\documentclass{article}
\usepackage{amsmath,amsfonts,amssymb,amsmath,amscd}
\usepackage[active]{srcltx}
\usepackage{graphicx,array}     

\newtheorem{theorem}{Theorem}[section]

\newtheorem{proposition}[theorem]{Proposition}

\newcommand{\dem}{\noindent{\it Proof:\ }}

\newcommand{\eproof}{\noindent\mbox{\framebox [0.6ex]{}}  \medskip}

\def \R{\mbox{${\mathbb R}$}} 
\def \s{\mbox{\rm{${\textbf{S}}$}}}

\def \z{\mbox{\rm{${\textbf{Z}}$}}}
\def \GL{\mbox{\rm{${\textbf{GL}}$}}}
\def \M{\mbox{${\mathcal{M}}$}}
\def \A{\mbox{${\mathcal{A}}$}}

\def \G{\mbox{${\mathcal{G}}$}}
\def \GG{\mbox{$\tilde{\mathcal{G}}$}}
\def \P{\mbox{${\mathcal{P}}$}}
\def \E{\mbox{${\mathcal{E}}$}}        

\def \PP{\mbox{${\vec{\mathcal{P}}}$}}
\def \QQ{\mbox{${\vec{\mathcal{Q}}}$}}
\def \EE{\mbox{${\vec{\mathcal{E}}}$}}
\def \FF{\mbox{${\vec{\mathcal{F}}}$}}

\newcommand{\adl}{Ad_L}
\newcommand{\adltk}{Ad_{L^t}^{\:k}}
\newcommand{\adlk}{Ad_{L}^{\:k}}

\DeclareMathOperator{\im}{Im}

\begin{document}

\title{{\bf Complete transversals of
reversible equivariant singularities of vector fields }}

\maketitle

\begin{center}
{\scshape Miriam Manoel} \\
 {\footnotesize {University of S\~ao Paulo, Mathematics Department}} \\
 {\footnotesize{S\~ao Carlos, SP - Brazil }} \\
 {\footnotesize{miriam@icmc.usp.br}} \\

{\scshape Iris de Oliveira Zeli} \\
{\footnotesize {Mathematics Department, IMECC, University of Campinas}} \\
{\footnotesize{Campinas, SP - Brazil}} \\
{\footnotesize{irisfalkoliv@ime.unicamp.br}} \\ 
\end{center}



\begin{abstract}
We use group representation theory to give algebraic formulae to compute complete transversals of singularities of vector fields, either in the nonsymmetric or in the reversible equivariant contexts. This computation produces normal forms directly, which are used sistematically in the local analysis of symmetric dynamics.
\end{abstract}



\section{Introduction}

In singularity theory there has been many results concerned with determining normal forms of map germs defined on different domains
under different equivalence relations.   Among a great number of papers in this direction, we cite for example the classical works by Bruce {\it et al.} \cite{bpw}, Gaffney and du Plessis \cite{gp}, Gaffney \cite{gaffney} and Wall \cite{wall,wall1}. On the classification of singularities  applied to bifurcation theory we mention Go\-lu\-bitsky {\it et al.} \cite{golub0,golub} and Melbourne \cite{melbourne,melbourne2}, these in the contexts with and without symmetries. In \cite{bkp} the authors present the so-called  {\it complete transversal method}, an important algebraic tool for the classification of finitely determined map germs.  In \cite{kirk} Kirk presents a programme in  {\it Maple}, called  Transversal, that implements
this method. \\

In dynamical systems, normal form theory of vector fields are obtained up to conjugacy and it has been an important tool for the study of local dynamics  around a singularity.  Some classical works are due to Poincar\'e \cite{poincare}, Birkhoff \cite{birkhoff}, Dulac \cite{dulac}, Belitskii \cite{belitskii}, Elphick {\it et al.} \cite{elphick} and Takens \cite{takens}. The method developed by Belitskii \cite{belitskii} consists of calculating the kernel of the homological operator associated with the adjoint $L^t$ of the linearization $L$ of the original vector field.  This calculation in turn  is associated with finding polynomial solutions of a PDE. Elphick {\it et al.} in \cite{elphick} give  an algebraic method for obtaining the normal form introducing an action of a group of symmetries $\s$, namely
\begin{equation}
\s=\overline{\left\{ e^{ sL^t}, s\in \R \right\}},
\label{eq:grupo-s}
\end{equation}
so that the polynomial nonlinear terms are equivariant under this action.
In \cite{bmz} we drice attention to formal normal forms of smooth vector fields in simultaneous presence  of symmetric and reversing symmetric transformations. The algebraic treatment shows advantage at once, since the set $\Gamma$ formed by such transformations  has a group structure. As a consequence, the vector field, called $\Gamma-$reversible-equivariant, has a well-determined general form that can be
 explicitly described by an algorithmic way, by means of results  given in \cite{abdm} and \cite{bm2}.
Purely reversible systems have been studied for a longer time, for example \cite{blm,champneys,jt1,jt2,lamb1,ltw}. In more recent years, reversible equivariant systems have become an object of study by many authors, see for example  \cite{abdm,bm2,bm1,bl,lr,lt2,mt}. In particular, in  \cite{bm1}
a relationship between purely equivariant  systems (without reversing symmetries) and a class of reversible equivariant
systems is established.
The normal form of a $\Gamma-$reversible-equivariant system  inherits the symmetries and reversing symmetries if the changes of coordinates are equivariant under the group $\Gamma$. Belitskii normal form has been used by many authors in different aspects; for example, in the analysis of occurrence of limit cycles or families of periodic orbits either in purely reversible vector fields or in reversible equivariant ones (see \cite{lt,lt2,mt,mereu}).
Motived by these works, in \cite{bmz} we have established an algebraic result related to those by Belitskii \cite{belitskii} and Elphick \cite{elphick} in the  reversible equivariant context using tools from invariant theory. In this process we have proved that the $\Gamma-$reversible equivariant normal form comes from the description of the reversible equivariant theory of the semidirect  product $\s \rtimes \Gamma$. After that recognition, we use results of \cite{abdm,bm2}  to produce the formal normal form of a reversible equivariant vector field by means of an alternative  algebraic method, without passing through a search for solutions of a PDE, which is the basis of Belitskii´s method.\\

In the present work we put together the approaches from singularity theory and  from dynamical systems theory in the study of normal forms. We show how the complete transversal method is closely related to the normal form method developed in \cite{bmz}. Let us stress that our intention here is not to apply the method for specific classifications. The goal is, instead, to explore this relation to recognize an algebraic alternative to compute complete transversals of singularities. Clearly the result is also  valid without symmetries. The idea is to introduce  two Lie groups of changes of coordinates, one for each context, and through these objects  we calculate the tangent spaces to the orbits of map germs under the appropriate actions of the groups. In the nonsymmetric context we recognize the complete transversal as being the space of polynomial map germs that commute with the group $\s$; in the reversible equivariant context, that turns out to be the space of polynomial map germs that are reversible equivariant under the action of  $\s \rtimes \Gamma$. \\

We have organized  this paper as follows. In Section \ref{section:preliminaries} we briefly present notation and basic concepts from singularity theory, from  reversible equivariant theory  and from normal form theory. In Section~\ref{section:algebraic transversal} we present the algebraic way to compute complete transversals.  According to the action of a specific group of equivalences, we characterize the tangent space to the orbit of a map germ (Proposition \ref{prop:tangent.space.G}),  and then we recognize the complete transversal (Theorem \ref{thm:transversal.for.G}). The main result appears in Subection~\ref{section:with.symmetry}, where we present the reversible equivariant version: we first consider a specific group of equivalences that preserves all symmetries and reversing symmetries of the original system, characterize the tangent space to orbit of a map germ (Proposition \ref{prop:tangent.space.GG}),  and  recognize the complete transversal in a purely algebraic way (Theorem \ref{thm:transversals.GG}).

\section{Preliminaries}\label{section:preliminaries}

In this section, we present the basic concepts from  singularity theory and reversible equivariant theory to present the main results.  We use the language of germs from singularity theory for the local study of $C^\infty$ applications around a singularity, which we assume to be the origin.

\subsection{Reversible equivariant theory}

Let $\E_V$ denote the space of all  smooth function germs $ V,0 \to \R$ and $\EE_{V}$ the space of all smooth map germs $V,0 \to V$. Let $\Gamma$ be a compact Lie group with a linear action  on a finite-dimensional real vector space $V$: $\Gamma \times V \to V, $ $(\gamma, x) \mapsto \gamma x$.

Consider a group homomorphism
\begin{equation}
\label{defisigma}
\sigma: \Gamma \to \z_2 =\{\pm 1 \},
\end{equation}
defining elements of $\Gamma$ as follows: if $\sigma(\gamma)= 1$ then $\gamma$ is a symmetry, if $\sigma(\gamma)=-1$, then $\gamma$ is a reversing symmetry. We denote by $\Gamma_+$ the subgroup of symmetries of $\Gamma$. If $\Gamma_+$ is nontrivial, then $\Gamma_+ = \ker \sigma$ is a
proper normal subgroup of $\Gamma$ of index $2$.

We recall that to a linear action of $\Gamma$ on $V$ there corresponds a representation $\rho$ of the group $\Gamma$ on $V$. In other words, there is a linear group homomorphism $\rho: \Gamma \to \GL(V)$, $\rho(\gamma)x = \gamma x$,  where $\GL(V)$ is the vector space of invertible linear maps $V \mapsto V$. The representation $\rho_\sigma: \Gamma \to \GL(V),$ $\rho_\sigma(\gamma) = \sigma(\gamma) \rho(\gamma)$, where $\sigma$ is the one-dimensional representation of $\Gamma$ given in \eqref{defisigma} is called the dual of $\rho$.

A smooth function germ $f: V,0 \to \R$ is called  $\Gamma-$invariant if
\begin{equation}
f(\rho(\gamma) x)=f(x), ~\forall \gamma \in \Gamma, \ x \in V,0.
\end{equation}

We denote by $\P_V(\Gamma)$ the ring of $\Gamma-$invariant polynomial function germs and by $\E_V(\Gamma)$ the ring of $\Gamma-$invariant smooth function germs.

A smooth map germ $g: V,0 \to V$ is called (purely) $\Gamma-$equivariant if
\begin{equation}
g(\rho(\gamma) x)=\rho(\gamma) g(x),~\forall \gamma \in \Gamma, \ x \in V,0.
\end{equation}
We denote by $\PP_V(\Gamma)$ the module of $\Gamma-$equivariant polynomial map germs and by $\EE_V(\Gamma)$ the module of $\Gamma-$equivariant smooth map germs.

A  smooth map germ $g: V,0 \to V$ is called $\Gamma-$reversible-equivariant if
\begin{equation}
g(\rho(\gamma) x)= \rho_\sigma(\gamma) g(x),~\forall \gamma \in \Gamma, \ x \in V,0.
\end{equation}

We denote by $\QQ_{V}(\Gamma)$ the module of $\Gamma-$reversible-equivariant polynomial map germs and by $\FF_V(\Gamma)$ the module of $\Gamma-$reversible-equivariant smooth map germs.

Since $\Gamma$ is compact, $\PP_{V}(\Gamma)$ and $\QQ_V(\Gamma)$ are finitely generated modules over $\P_V(\Gamma),$ which in turn is  finitely generated ring (see \cite{abdm,golub}). If $\sigma$ is trivial, then  $\PP_V(\Gamma)$ and $\QQ_V(\Gamma)$ coincide. In \cite{abdm}, the authors present an algorithm that produces a generating set of $\QQ_V(\Gamma)$ over $\P_{V}(\Gamma)$.  A result in \cite{bm2} provides a simple way to compute a set of generators of  $\P_V(\Gamma)$ from the knowledge of generators of $\P_V(\Gamma_+)$.

Notice that $\PP_V(\Gamma)$ and $\QQ_V(\Gamma)$ are graduate modules,
\begin{equation}
\PP_V(\Gamma) = \bigoplus_{k \geq 0} \PP_V^k(\Gamma) \quad \text{and} \quad \QQ_V(\Gamma) = \bigoplus_{k \geq 0} \QQ_V^k(\Gamma),
\end{equation}
for $\PP_V^k(\Gamma) = \PP_V(\Gamma) \cap \PP_V^k~$ and $~\QQ_V^k(\Gamma) = \QQ_V(\Gamma) \cap \PP_V^k$, where $\PP_V^k$ is the space of homogeneous polynomials of degree $k$ defined on $V$, $k \geq 0$.

\subsection{Belitskii-Elphick normal form}

We begin with a brief introduction to the Belitskii method. For $h \in \EE_V$, consider the ODE
\begin{equation}
\label{eq:sistemainicial}
\dot{x}= h(x),~ x\in V,0.
\end{equation}

The interest of the theory is local, around a singular point which we assume to be the origin, so $h(0)=0$. The  method to calculate the normal form consists of  successive changes of coordinates in the domain that are perturbations of the identity, i.e.,  of type $x=\xi(y)= y+ \xi_k(y)$, $k \geq 2$,  where $I$ is the germ of the identity and $\xi_k \in \PP_V^k$. In the new variables, the system is
$$
\dot{y} = g(y),~y\in V,0,
$$
where
\begin{equation}
\label{eq:campo-conjugado}
g(y)=  (d\xi)_x^{-1} h (\xi(y)),
\end{equation}

Note that for each $x$ we have
\begin{equation}
\label{eq:inversa.mudanca}
(d\xi)_x^{-1}=(I + (d\xi_k)_x)^{-1} = I - (d\xi_k)_x + \varphi \left( (d\xi_k)_x \right), ~ k\geq 2,
\end{equation}
where $\varphi((d\xi_k)_x)$ is polynomial and has degree greater or equal to $2(k-1)$.

The aim is to annihilate  as many terms of degree $k$ as possible in the original vector field, obtaining a conjugate vector field  written in a simpler and more convenient form.  The method developed by Belitskii reduces this problem to computing $\ker \adltk$ where  $\adlk: \PP_V^k \to \PP_V^k$ is the so-called homological operator and is defined by
\begin{equation}
\label{eq:operatorhomo}
\adlk(p)(x) = (d p)_{x}Lx-Lp(x),~~x \in V,0,
\end{equation}
being $L^t$ the adjoint of the linearization  $L$. We refer to \cite{golub} for the details.

In \cite{elphick}, Elphick {\it et al.} give an alternative algebraic method to obtain the  normal form developed by Belitskii, which consists
 of computing nonlinear terms that are equivariant under the action of the group $\s$ given by
\begin{equation}
\s=\overline{\left\{ e^{ sL^t}, s\in \R \right\}}.
\label{eq:group.S}
\end{equation}
In fact, the authors show that for each $k \geq 2$, $\ker \adltk =\PP_V^k(\s)$ and, therefore,
\begin{equation}
\label{eq:elphick}
\PP_V^{k}= \PP_V^k(\s) \oplus \adlk(\PP_V^k).
\end{equation}
From that, we show in \cite[Theorem 4.1]{bmz}, that if the vector field $h$ is $\Gamma-$reversible-equivariant, with $L= (dh)_0$, then  for each $k \geq 2$ we have
\begin{equation}
\label{eq:decomposition.nf.rever.equiv}
\QQ_V^{k}(\Gamma)= \QQ_V^{k}(\s \rtimes \Gamma) \oplus \adl^{k}  ( \PP_V^{k}(\Gamma)),
\end{equation}
where the semidirect product is induced from the homomorphism $\mu: \Gamma \to Aut(\s)$ given by
$$
\mu(\gamma)(e^{sL^t})= e^{\sigma(\gamma)L^t}.
$$
Hence, the normal form deduction reduces to the computation of a basis for the vector space $\QQ_V^{k}(\s \rtimes \Gamma)$ for each $k \geq 2$.  In practice, we obtain the general form of elements in  $\QQ_V(\s \rtimes \Gamma)$ and, once this module is graduated, we easily extract from that a basis for
 $\QQ_V^{k}(\s \rtimes \Gamma)$. The main tools we use to obtain this general form  are \cite[Algorithm 3.7]{abdm} and \cite[Theorem 3.2]{bm2} which hold in particular if the group is compact. There are many cases for which the group $\s$ fails to be compact; nevertheless,  these tools can still be used as long as the ring $\P_V(\s)$ and the module  $\PP_V(\s)$ are finitely generated.


\section{The algebraic alternative for complete transversals} \label{section:algebraic transversal}

\subsection{The nonsymmetric context}

Recall that $\E_V$ is a commutative local ring with  maximal ideal given by
\begin{equation}
\label{eq:maximal.ideal}
\M = \left\{ f\in \E_V, f(0)=0 \right\}.
\end{equation}
For each $k$, the power  $\M^{k}$ is the ideal of function germs in $\EE_V$ whose Taylor polynomials up to degree $k$ vanish at the origin.
We shall also consider the vector space $J^k$ formed by all $k-$jets $j^kh$ of elements $h \in \M \EE_V$.

 We denote by  $\G$ the group of  germs of diffeomorphisms $\xi: V,0 \to V$, $~\xi = I + \tilde{\xi}$, where  $I$ is the germ of identity and  $\tilde{\xi} \in  \bigoplus_{l \geq 2} \PP_V^l$, which is  a  geometric subgroup in the sense of Damon \cite{damon}. We  consider the action of $\G$ on $\M \EE_V$ given by
\begin{equation}
\label{eq:action.G}
\left(\xi \cdot h  \right)(x) = (d\xi)_{\xi(x)}^{-1} h  \left( \xi (x)\right),~\xi \in \G,~h\in \M\EE_V.
\end{equation}

For each $k\geq 2$ we introduce the group  $J^k\G \ = \ \left\{j^k \xi,~\xi \in \G  \right\}$, which is a Lie group with an action on $J^k$  induced by  \eqref{eq:action.G}, namely
$$
j^k\xi \cdot (j^kh)  (x) = j^k( \xi \cdot h)(x),~~\xi \in \G,~~h\in \M \EE_V.
$$
For this action, we consider the orbit of $h$ and define  the tangent space  $T \G \cdot h$  to this orbit by  the set  of elements of the
form
\begin{equation}
\frac{d}{dt} \phi(x,t)_{|t=0},
\end{equation}
for the one-parameter family $\phi(\cdot,t)$ where $\phi(x,t)=  (d\xi)^{-1}_{(x,t)} h \left(\xi(x,t) \right)$ with  $\xi(x,0)=x$.

The complete transversal method by Bruce {\it et al.} \cite{bkp} is a tool for the classification of singularities that is performed on each degree level
in the Taylor expansion of the germ to be studied. The main idea is to classify, on each step, $k-$jets on $J^k$,  once $J^k$  is isomorphic to a quotient of $\E_V-$modules $\M \EE_V / \M^{k+1}\EE_V ~$. In the proposition below we transcribe the result (\cite[Proposition 2.2]{bkp}):
\begin{proposition} \label{thm:transversals} For $k \geq 1$, let $h$ be a $k-$jet in the jet space $J^k$. If $W$ is a vector subspace of $\PP_V^{k+1}$  such that
\begin{equation}
\label{eq:complete.transversal}
\M^{k+1}\EE_V \subset W + T\G\cdot h +  \M^{k+2} \EE_V,
\end{equation}
then every $k+1-$jet $g$ with $j^kg= h$ is in the same $J^{k+1}\G-$orbit as some $(k+1)-$jet of the form $h + \omega$, for some $\omega \in W$.
\end{proposition}

The vector subspace $W$ is the so-called complete transversal. In principle, the computation of $W$ requires the knowledge of   $T\G\cdot h$
modulo $\M^{k+2} \EE_V$. Now, in an investigation of this result, we have noticed the presence of a linear operator  resembling
the homological operator given in \eqref{eq:operatorhomo}. This has lead us to obtain an alternative way to compute complete transversals
through an algebraic approach. The rest of the present work is devoted to present that approach.

Let us  introduce the linear operator $Ad_{h}: \EE_V \to \EE_V$ which is a generalization of the homological operator and is defined by
\begin{equation}
\label{eq:linear.operator}
Ad_{h}(\xi)(x)= (d\xi)_x h(x) - (dh)_x \xi(x).
\end{equation}

Now, we consider the restriction $Ad_{h}^{\:k}= {Ad_{h}}|_{\PP_V^k}$. Writing $h$ as  $h = L + \tilde{h}$ \  with $L= (dh)_0$ and $\tilde{h} \in \M^{2} \EE_V$, \ it follows from the linearity of $Ad_{h}$ that
\begin{equation}
Ad_h^{\:k}(\xi_k)= \adlk(\xi_k) + Ad_{\tilde{h}}^{\:k}(\xi_k),~~\xi_k \in \PP_V^k.
\end{equation}
From that  we can now characterize the tangent space $T\G \cdot h$:

\begin{proposition} \label{prop:tangent.space.G} The tangent space to the orbit of $h \in \M \EE_V$ is given by
$$
T\G \cdot h =  \left\{ Ad_{h}(\tilde{\xi}) + {\varphi}(-(d\tilde{\xi})_x) h, ~ ~\tilde{\xi} \in \bigoplus^{l \geq k}\PP_V^l, ~\varphi((d\tilde{\xi})_x) ~\text{as in}~ \eqref{eq:inversa.mudanca},~k \geq 2 \right\} .
$$
\end{proposition}

\dem Let $\xi(\cdot,t)$ be a family on $\G$, $\xi(x,t)=x + \tilde{\xi}(x,t)$, with $\xi (x,0)=x$, and let
$$
\phi(x,t)= (d\xi)^{-1}_{\xi(x,t)}h\left( \xi (x,t)\right).
$$
We have
$$
\frac{d}{dt}{\phi}(x,0)= \left( -\frac{d}{dt}{(d \tilde{\xi})}_{x} + \varphi \left( \frac{d}{dt} (d\tilde{\xi})_x \right) \right)h(x) + (dh)_{x} \frac{d}{dt}\tilde{\xi}(x,0),
$$
with $\varphi$  given by
\begin{equation}
\label{eq:varphi}
(d\xi)^{-1}_{\xi(x,t)}= I - (d \tilde{\xi})_{\xi(x,t)}+ \varphi((d\tilde{\xi})_{(x,t)}).
\end{equation}
Rewriting
\begin{equation}
\label{eq:varphi.identified}
\displaystyle\frac{d}{dt}{(d \tilde{\xi})}_{x}\equiv(d \tilde{\xi})_x ~, \quad  \varphi \left( \frac{d}{dt}(d \tilde{\xi})_x \right) \equiv \varphi((d\tilde{\xi})_x)  \quad  \text{and} \quad  \displaystyle\frac{d}{dt}\tilde{\xi}(x,0)\equiv \tilde{\xi}(x),
\end{equation}
the result follows immediately.
\eproof

The theorem below is now a direct consequence of Proposition~ \ref{prop:tangent.space.G}:
\begin{theorem} \label{thm:transversal.for.G} For $k \geq 1$ let $h \in J^k$. Consider the vector subspace $\PP_V^{k+1}(\s)$ of  $\PP_V^{k+1}$,
with $\s$ defined in \eqref{eq:group.S} associated with $L=(dh)_0$. Then,
$$
\M^{k+1}\EE_V \subset \PP_V^{k+1}(\s) + T \G \cdot h +  \M^{k+2}\EE_V.
$$

\end{theorem}
\dem Let $g\in \M^{k+1}\EE_V$. From the  decomposition \eqref{eq:elphick}, for each degree-$k$ term  $g_{k+1}$ in the Taylor expansion of $g$ we have
$$
g_{k+1} = q_{k+1}~+~p_{k+1},
$$
 with  $q_{k+1} \in \PP_V^{k+1}(\s)~$   and  $~p_{k+1} \in \im \adl^{\:k+1}$.
Then,  $p_{k+1} = \adl^{k+1}(\xi_{k+1})$ for some  $\xi_{k+1} \in \PP_V^{\:k+1}$. Consider $\varphi(-(d\xi_{k+1})_x)$ as in  \eqref{eq:inversa.mudanca}. We write  $h=L +\tilde{h}$, with  $L=(dh)_0$ and $\tilde{h} \in \M^{2} \EE_V$, to obtain
\begin{equation} \nonumber
g_{k+1} =  q_{k+1} +  Ad_{h}(\xi_{k+1}) + \varphi(-(d\xi_{k+1})_x)  h  - \left( Ad_{\tilde{h}}(\xi_{k+1}) + \varphi(-(d\xi_{k+1})_x)  h \right). \nonumber
\end{equation}

By Proposition \ref{prop:tangent.space.G}, $~Ad_{h}(\xi_{k+1})$ + $\varphi\left( -(d\xi_{k+1})_x \right)  h \in T\G \cdot h$. Furthermore, from the definition of the linear operator and $\tilde{h}$ it follows that
$$ Ad_{\tilde{h}}(\xi_{k+1}) + \varphi(-(d\xi_{k+1})_x) h \in  \M^{k+2} \EE_V.$$
\eproof

We remark that the choice of a vector subspace $W$ satisfying \eqref{eq:complete.transversal} is obviously not unique; however, by decomposition \eqref{eq:elphick}, $\PP_V^k(\s)$ is such a  space which has the smallest dimension.


\subsection{The reversible equivariant context} \label{section:with.symmetry}

Let $\Gamma$ be a compact Lie group and consider the homomorphism $\sigma$ defined in \eqref{defisigma}. Here we extend the results of the previous subsection to the  $\Gamma-$reversible-equivariant context. In particular, if  $\sigma$ is trivial then the result reduces to the (purely) $\Gamma-$equivariant context. \\

Let us denote by $\GG$ the group formed by formal changes of coordinates of type $\xi: V,0 \to V$, $\xi = I + \tilde{\xi}$, where  $I$ is the germ of the identity and  $\tilde{\xi} \in \bigoplus_{l \geq 2} \PP_V^l(\Gamma)$.
 $\GG$ is a subgroup of $\G$ and the action of $\GG$ on $\FF_V(\Gamma)$ is defined as in \eqref{eq:action.G}.

Our space of germs is now $\FF_V(\Gamma)$.  Let us denote by $J^k (\Gamma_\sigma)$ the space of $\Gamma-$reversible-equivariant $k-$jets and,
for each $k \geq 1$, we denote by  ${\FF_V}_{k+1}(\Gamma)$ the space $\M^{k+1}\EE_V  \cap \FF_V(\Gamma)$.  Also, for each $k \geq 1$,
let $J^k\GG$ denote the  group of $k-$jets $j^k \xi$ of elements $\xi \in \GG$. Consider now the action of $J^k\GG$ on  $J^k(\Gamma_\sigma)$ induced by \eqref{eq:action.G}:
$$
j^k \xi\cdot \left( j^kh\right) (x) = j^k( \xi \cdot h)(x),~~\xi \in \GG,~h\in \FF_V(\Gamma),~h(0)=0.
$$

Castro and du Plessis have stated in \cite{cp}  the equivariant version of Proposition~\ref{thm:transversals}.
We have verified that the reversible equivariant version  is direct if we consider the group $\GG$. The statement of the theorem  is:

\begin{theorem} \label{thm:transversals.GG} For $k \geq 1$ let $h$ be a $k-$jet in the jet space $J^k(\Gamma_\sigma)$. If \  $W$ is a vector subspace of $\QQ_V^{k+1}(\Gamma)$  such that
\begin{equation}
\label{eq:thm:transversals.GG}
{\FF_V}_{k+1}(\Gamma)  \subset W + T\GG\cdot h +  {\FF_V}_{k+2}(\Gamma),
\end{equation}
then every $\Gamma-$reversible-equivariant $k+1-$jet $g$ with $j^kg= h$ is in the same $J^{k+1}\GG-$orbit as some $(k+1)-$jet of the form $h + \omega$, for some $\omega \in W$.
\end{theorem}
As in the previous subsection, our aim here is to determine a subspace $W$ satisfying \eqref{eq:thm:transversals.GG}. For that, we first
characterize the tangent space $T\GG \cdot h$ for $h\in \FF_V(\Gamma)$, $h(0)=0$ through   the linear operator defined in \eqref{eq:linear.operator}:

\begin{proposition} \label{prop:tangent.space.GG} For $h \in \FF_V(\Gamma)$ with $h(0)=0$, the tangent space to the orbit of $h$ is given by
$$
T\GG \cdot h = \left\{ Ad_{h}(\tilde{\xi}) + \varphi((d\tilde{\xi})_x) h,~~\tilde{\xi} \in \bigoplus^{l\geq k}\PP_V^l(\Gamma), ~ \varphi((d\tilde{\xi})_x) ~\text{as in}~ \eqref{eq:inversa.mudanca},~k\geq2 \right\}.
$$
\end{proposition}

The proof of this proposition follows the steps of the proof of Proposition \ref{prop:tangent.space.G}, accompanied with the $\Gamma$-equivariance. \\

The result below reveals the alternative complete transversal for the reversible equivariants:

\begin{theorem} \label{thm:transversal.for.GG} For $k \geq 1$, let $h \in J^k (\Gamma_\sigma)$, $L=(dh)_0$. Consider the group $\s$ given in \eqref{eq:group.S} associated with $L$.   Then,
$$
{\FF_V}_{k+1}(\Gamma)  \subset \QQ_V^{k+1}(\s\rtimes\Gamma) +  T \GG \cdot h + {\FF_V}_{k+2}(\Gamma).
$$
\end{theorem}

\dem Use the decomposition \eqref{eq:decomposition.nf.rever.equiv} and follow the steps of the proof of Theorem \ref{thm:transversal.for.G}.
\eproof

As in the context without nontrivial symmetries, $\QQ_V^k(\s \rtimes \Gamma)$  is  a complete transversal  of smallest dimension that satisfies \eqref{eq:thm:transversals.GG}.\\

{\it Acknowledgments}: We wish to thank M. Escudeiro for a helpful suggestion.


\end{document}